\documentclass [11pt,a4paper]{article}
\usepackage{amssymb}
\def\F{\mathbb{F}}

\def\Z{\mathbb{Z}}

\usepackage{color,xcolor}
\newtheorem{theorem}{Theorem}
\newtheorem{lemma}{Lemma}

\newcommand{\quash}[1]{}

\begin{document}

\title{Trace representation and linear complexity of binary
sequences derived from Fermat quotients}

\author{Zhixiong Chen\\
Department of Mathematics, Putian University, \\ Putian, Fujian
351100, P.R. China\\
ptczx@126.com\\
}
\maketitle

\begin{abstract}
We describe the trace representations of two families of binary sequences
derived from Fermat quotients modulo an odd prime $p$ (one is the binary threshold sequences, the other is the Legendre-Fermat quotient sequences) via determining the defining pairs of all binary characteristic sequences of cosets, which coincide with the sets of pre-images modulo $p^2$ of each fixed value of Fermat quotients. From the  defining pairs, we can obtain an earlier result of linear complexity for the binary threshold sequences and a new result of linear complexity for the Legendre-Fermat quotient sequences  under the assumption of $2^{p-1}\not\equiv 1 \bmod {p^2}$.
\end{abstract}



\section{Introduction}\label{intro}
For an odd prime $p$ and an integer $u$ with $\gcd(u,p)=1$, the {\it
Fermat quotient\/} $q_p(u)$ modulo~$p$ is defined as the unique
integer with
$$
q_p(u) \equiv \frac{u^{p-1} -1}{p} \pmod p, \qquad 0 \le q_p(u) \le
p-1,
$$
and we also define
$$
q_p(kp) = 0, \qquad k \in \mathbb{Z}.
$$
Many number theoretic and cryptographic questions as well as
measures of pseudorandomness have been studied for Fermat quotients
and their generalizations
\cite{ADS,AW,BFKS,C,CD,CG,COW,CW2,CW,CW3,DCH,DKC,EM,GW,OS,Sha,Shk,S,S2010,S2011,S2011b,SW}.

In this paper, we still concentrate on certain binary sequences defined from Fermat quotients in the references.
The first one is the binary threshold sequence $(e_u)$ studied in \cite{CD,COW} by defining
\begin{equation}\label{binarythreshold}
e_u=\left\{
\begin{array}{ll}
0, & \mathrm{if}\,\ 0\leq q_p(u)/p< \frac{1}{2},\\
1, & \mathrm{if}\,\ \frac{1}{2}\leq q_p(u)/p< 1,
\end{array}
\right. \quad  u \ge 0.
\end{equation}
The second one, by combining $q_p(u)$ with the Legendre symbol
$\left( \frac{\cdot}{p} \right)$, is defined in \cite{GW} by
\begin{equation}\label{binarylegendre}
f_u=\left\{
\begin{array}{ll}
0, & \mathrm{if}\,\ \left(\frac{q_p(u)}{p}\right)=1\,\ \mathrm{or}\,\ q_p(u)=0, \\
1, & \mathrm{otherwise},
\end{array}
\right. \quad u \ge 0.
\end{equation}
(In fact, in \cite{GW} $\chi$, a fixed multiplicative
character modulo $p$ of order $m>1$, is applied to defining $m$-ary sequences $(h_u)$ of
discrete logarithms modulo a divisor $m$ of $p-1$ by
$$
  \exp(2\pi i h_u/m)=
     \chi(q_p(u)),~ 0\le h_u<m \quad \mbox{if } q_p(u)\not\equiv 0\bmod p
$$
and $h_u=0$ otherwise. When $m=2$, we have $h_u=f_u$ for all $u\ge 0$.)  We note that $(e_u)$ and $(f_u)$ are $p^2$-periodic since
\begin{equation}\label{fqprop}
q_p(u+kp) \equiv q_p(u) - ku^{-1} \pmod p, ~~\gcd(u,p)=1,
\end{equation}
see, e.g. \cite{OS}.

The authors of \cite{COW,GW} investigated measures of pseudorandomness as well as linear complexity profile of
$(e_u)$ and $(h_u)$ (of course including $(f_u)$) via certain character sums over Fermat quotients, respectively. The authors of \cite{CD}  determined the \emph{linear complexity} (see
 Section \ref{prelim} for the definition) of $(e_u)$ under the condition  $2^{p-1}\not\equiv 1 \pmod {p^2}$ by using the theory of cyclotomy, since $q_p: \Z_{p^2}^*\rightarrow \Z_p$ is a group homomorphism by the fact, see e.g. \cite{OS}, that
\begin{equation}\label{fqpropadd}
q_p(uv) \equiv q_p(u) + q_p(v) \pmod p, ~~\gcd(uv,p)=1.
\end{equation}
However, the linear complexity of $(f_u)$ is still open. (Note that part work has been done in \cite{CHD} if $2$ is a primitive element
modulo $p^2$.) One of our main aims in this paper  is to solve this problem.
Our second main aim is to investigate the trace representations of $(e_u)$ and $(f_u)$.  For our purpose, we need to describe $(e_u)$  and $(f_u)$ in an  equivalent  way.

Define
$$
D_l=\{u: 0\le u< p^2,~ \gcd(u,p)=1,~ q_p(u)=l\}
$$
for $l=0,1,\ldots,p-1$. Indeed, if $g$ is a (fixed) primitive root modulo $p^2$, we have by (\ref{fqpropadd})
$$
D_0=\{g^{kp} \bmod {p^2}: 0\le k<p\}
$$
and for $l=0,1,\ldots,p-1$, there exists an integer $0\le l_0<p$ such that $D_l=g^{l_0}D_0$. Hence each $D_l$ has the cardinality $|D_l|=p-1$. We remark that
$l_0=l$ if $q_p(g)=1$. Naturally $D_0,D_1, \ldots, D_{p-1}$ form a partition of $\mathbb{Z}_{p^2}^*$.

Let $P=\{kp: 0\le k< p\}$,
one can equivalently define
$(e_u)$  and $(f_u)$ respectively by resetting the terms in the first period and repeating periodically
$$
e_u=\left\{
\begin{array}{ll}
0, & \mathrm{if}\,\ u\in D_0 \cup \cdots \cup D_{(p-1)/2} \cup P,\\
1, & \mathrm{if}\,\ u\in D_{(p+1)/2} \cup \cdots \cup D_{p-1},
\end{array}
\right. \quad 0\le u < p^2
$$
and
$$
f_u=\left\{
\begin{array}{ll}
0, & \mathrm{if}\,\ u\in \cup_{l\in Q} D_l  \cup D_0\cup P,\\
1, & \mathrm{if}\,\ u\in \cup_{l\in N} D_{l},
\end{array}
\right. \quad 0\le u < p^2,
$$
where $Q$ is the set of quadratic residues modulo $p$ and $N$ the set of quadratic non-residues modulo $p$. We note that
the cardinality $|Q|=|N|=(p-1)/2$.

Trace function is extensively applied to
producing pseudorandom sequences efficiently and analyzing their pseudorandom properties.
A well-known example is the $m$-sequences \cite{GG,LN}. The trace representations of many famous sequences, such as Legendre and Jacobi sequences and their generalizations, have been studied in the literature, see  \cite{DGS,DGS09,DGSY,DYX,HKN,KS,KSG,NLCSY}. In this paper, we will represent $(e_u)$ and $(f_u)$ as a sum of trace functions via determining the \emph{defining pair} (see
Section \ref{prelim} for the definition) of the binary sequences $(s^{(l)}_u)$ defined by
$$
s^{(l)}_u=\left\{
\begin{array}{ll}
1, & \mathrm{if}\,\ u \bmod {p^2}\in D_l,\\
0, & \mathrm{otherwise},
\end{array}
\right. \quad  u \ge 0,~~ l=0,1,\ldots,p-1,
$$
which is called the \emph{binary characteristic sequence} with respect to the coset $D_l$. As an application of  defining pairs,
we obtain  an earlier known result of linear complexity  of $(e_u)$ proved in \cite{CD} and a new result of linear complexity of $(f_u)$ under the assumption of $2^{p-1}\not\equiv 1 \pmod {p^2}$.

\section{Preliminaries}\label{prelim}

Let $\mathbb{F}_2=\{0,1\}$ be the binary field and $\overline{\F}_2$ the algebraic closure of $\mathbb{F}_2$. For  a binary sequence $(s_u)$ over $\F_2$ of odd period $T$, there exists a primitive $T$-th root
$\beta\in\overline{\F}_2$ of unity and a polynomial $g(x)= \sum\limits_{0\le i<T}\rho_ix^i \in \overline{\F}_2[x]$ such that
$$
s_u=g(\beta^u), ~~ u\ge 0,
$$
see \cite[Theorem 6.8.2]{J}, we call the pair  $(g(x),\beta)$ a \emph{defining pair} of $(s_u)$ and $g(x)$ the \emph{defining polynomial} of $(s_u)$ corresponding to $\beta$  \cite{DGS,DGS09,DGSY}. Note that for a given $\beta$, $g(x)$ is uniquely determined up to modulo $x^{T}-1$ \cite[Lemma 2]{DGSY}. The relation between $(s_u)$
 and  $\{\rho_i: 0\le i<T\}$ is given as
$$
s_u=\sum\limits_{0\le i<T}\rho_i\beta^{iu} \Longleftrightarrow \rho_i=\sum\limits_{0\le u<T}s_u\beta^{-iu}.
$$
The righthand side is referred to as the \emph{discrete Fourier transform} of $(s_u)$ \cite[Ch. 6]{GG}.

 We recall that the
\emph{linear complexity} $L((s_u))$  is the least order $L$ of a linear
recurrence relation over $\mathbb{F}_2$
$$
s_{u+L} = c_{L-1}s_{u+L-1} +\cdots +c_1s_{u+1}+ c_0s_u\quad
\mathrm{for}\,\ u \geq 0,
$$
which is satisfied by $(s_u)$ and where $c_0=1, c_1, \ldots,
c_{L-1}\in \mathbb{F}_2$. The linear complexity  of sequences plays an important role in stream cipher.
For a sequence to be cryptographically strong, its linear complexity
should be large and at least a half of the period  according to  the Berlekamp-Massey
algorithm \cite{Massey}. From \cite{B} or \cite[Theorem 6.3]{GG}, the linear complexity of $(s_u)$ is determined by
\begin{equation}\label{lc-rho}
L((s_u))=|\{i: \rho_i\neq 0, 0\le i<T\}|,
\end{equation}
i.e., the linear complexity of $(s_u)$ equals the number of nonzero coefficients of $g(x) \bmod x^{T}-1$, which is also called the \emph{Hamming weight} of $g(x)$.

However it is not easy to determine the linear complexity via (\ref{lc-rho})
by re-constructing $g(x)$ for a sequence. Fortunately one can determine the linear complexity in another way.
Let
$$
S(x)=s_0+s_1x+s_2x^2+\cdots+s_{T-1}x^{T-1}\in \mathbb{F}_2[x],
$$
which is called the \emph{generating polynomial} of $(s_u)$. Then the linear
complexity of $(s_u)$ is computed by
\begin{equation}\label{licom}
  L((s_u)) =T-\deg\left(\mathrm{gcd}(x^T-1,
  ~S(x))\right),
\end{equation}
see, e.g.  \cite{LN} for details.\\

First we present some auxiliary statements. Let $g$ be a (fixed) primitive root modulo $p^2$ and $q_p(g)=\delta$ for some $1\le \delta<p$ in the context. For the sake of convenience, we adjust the subscript of cosets $D_l$ defined in Section 1 in the following way:
\begin{equation}\label{Dj}
D_{j\delta}=g^jD_0=\{g^{kp+j} \bmod {p^2}: 0\le k<p\},~~ 0\le j<p,
\end{equation}
here and hereafter the subscript of $D$ is performed modulo $p$, i.e., $D_{l+p}=D_l$ for all $l\ge 0$.

\begin{lemma}\label{lemma1}\cite{CD}
For any $0\leq l< p$, if $u\bmod {p^2}\in D_{l'}$ for some $0\leq l'< p$, we
have
$$
u D_{l} =\{uv \bmod {p^2} : v\in D_{l} \} = D_{l+l'}.
$$
\end{lemma}

Define   $$ D_l(x)= \sum\limits_{u\in D_l}x^u \in \mathbb{F}_2[x]$$
for $0\leq l < p$.

\begin{lemma}\label{lemma-add}
Let $\beta \in \overline{\mathbb{F}}_{2}$ be  a  primitive $p^2$-th root of
unity.

(i). \cite{CD} For all $n\in \mathbb{Z}_{p^2}^*$, we have
$$\sum\limits_{l=0}^{p-1}D_l(\beta^n)=0.$$

(ii).  For $0\leq l < p$, we have
$$
D_l(\beta^{kp})=\left\{
\begin{array}{ll}
0, & \mathrm{if}~ k\equiv 0 \pmod p,\\
1, & \mathrm{otherwise}.
\end{array}
\right.
$$
\end{lemma}
Proof.  We only prove (ii). By  (\ref{Dj}), we have
$$
\{u\pmod p: u\in D_l\}=\Z_p^*, ~~~0\leq l < p.
$$
Let $\theta=\beta^p$, which is  a  primitive $p$-th root of
unity.  For $0\leq l < p$, we derive
$$
 D_l(\beta^{kp})=  D_l(\theta^{k})=\sum\limits_{u\in D_l}\theta^{ku}=\sum\limits_{j\in \Z_p^*}\theta^{kj},
 $$
which deduces the desired result for different $k$ modulo $p$. We note that the calculations are performed in
finite fields with characteristic two. ~\hfill $\square$

Motivated by \cite{DGSY} we describe a method for determining the defining pair of binary characteristic sequences with respect to a single coset in the following technical lemma, which will be used to show the main theorem.

\begin{lemma}\label{single-def-pair}
Let $q_p(g)=\delta$ for a (fixed) primitive root $g$ modulo $p^2$ and $\beta \in \overline{\mathbb{F}}_{2}$ be a  primitive $p^2$-th root of unity. Then for $0\le i<p$, the defining pair of $(s^{(i\delta)}_u)$ defined by
$$
s^{(i\delta)}_u=\left\{
\begin{array}{ll}
1, & \mathrm{if}\,\ u\bmod {p^2} \in D_{i\delta},\\
0, & \mathrm{otherwise},
\end{array}
\right. \quad  u \ge 0,
$$
is $(G_{i\delta}(x),\beta)$ with
$$
G_{i\delta}(x)=\sum\limits_{j=1}^{p-1}x^{jp}+\sum\limits_{j=0}^{p-1} D_{(i+j)\delta}(\beta)D_{j\delta}(x).
$$
The numbers $i\delta, j\delta, (i+j)\delta$ in the super- or sub-script are reduced modulo $p$.
\end{lemma}
Proof.
Define $p$-tuples
$$
C_i=(D_{i\delta}(\beta),D_{(i+1)\delta}(\beta),\ldots,D_{(i+p-1)\delta}(\beta)), ~~ i=0,1,\ldots,p-1.
$$
Using (\ref{Dj}) we calculate for fixed $0\le i,j<p$
\begin{eqnarray*}
C_iC^{T}_j& = & D_{i\delta}(\beta)D_{j\delta}(\beta)+D_{(i+1)\delta}(\beta)D_{(j+1)\delta}(\beta)+\ldots+D_{(i+p-1)\delta}(\beta)D_{(j+p-1)\delta}(\beta)\\
          & = & \sum\limits_{k=0}^{p-1}~\sum\limits_{u\in D_0}\beta^{ug^{i+k}}~\sum\limits_{v\in D_0}\beta^{vg^{j+k}}\\
          & = & \sum\limits_{k=0}^{p-1}~\sum\limits_{u\in D_0}\beta^{ug^{i+k}}~\sum\limits_{w\in D_0}\beta^{uwg^{j+k}}~~ ~~(\mathrm{we~ use~ } v=uw)\\
          & = & \sum\limits_{k=0}^{p-1}~\sum\limits_{u\in D_0}~\sum\limits_{w\in D_0}\beta^{ug^{j+k}(g^{i-j}+w)}\\
          & = & \sum\limits_{w\in D_0}~\sum\limits_{z\in \Z_{p^2}^*}\gamma_w^z ~~(\mathrm{we~ use~ } z=ug^{j+k}, \gamma_w=\beta^{g^{i-j}+w}).
\end{eqnarray*}
For those $w\in D_0$ with $\gcd(g^{i-j}+w, p)=1$, we see that  $\gamma_w$ is  a  primitive $p^2$-th root of
unity and hence we have
$$
\sum\limits_{z\in \Z_{p^2}^*}\gamma_w^z=\sum\limits_{l=0}^{p-1}D_l(\gamma_w)=0
$$
by Lemma \ref{lemma-add}(i).

For $w\in D_0$ with $\gcd(g^{i-j}+w, p)=p$, we suppose $g^{i-j}+w\equiv l_0p \pmod {p^2}$ for some integer $0\le l_0<p$ and use (\ref{fqprop}) to get
$$
0\equiv q_p(w)\equiv q_p(-g^{i-j}+l_0p)\equiv q_p(-g^{i-j})-l_0(-g^{i-j})^{-1} \pmod p,
$$
which implies that there is exactly one $w$ satisfying such condition. In particular, if $i=j$ we see that $l_0=0$ and hence
$\gamma_w=1$, in this case we have
$$
\sum\limits_{z\in \Z_{p^2}^*}\gamma_w^z=p(p-1)=0.
$$
However if  $i\neq j$, we see that $\gamma_w$ is a  primitive $p$-th root of
unity and by Lemma \ref{lemma-add}(ii) we derive
$$
\sum\limits_{z\in \Z_{p^2}^*}\gamma_w^z=\sum\limits_{l=0}^{p-1}D_l(\gamma_w)=p=1.
$$

Putting everything together, we get
$$
C_iC^{T}_j=\left\{
\begin{array}{ll}
0, & \mathrm{if}\,\ i=j,\\
1, & \mathrm{otherwise}.
\end{array}
\right.
$$
Now we verify that $(G_{i\delta}(x),\beta)$ is the defining pair of $(s^{(i\delta)}_u)$. For $u\equiv 0 \pmod {p^2}$ we have
$$
G_{i\delta}(\beta^u)=G_{i\delta}(1)=p-1+(p-1)\sum\limits_{j=0}^{p-1} D_{(i+j)\delta}(\beta)=0=s^{(i\delta)}_u.
$$
For $u\equiv lp \pmod {p^2}$ with $1\le l<p$, we have by Lemma \ref{lemma-add}
\begin{eqnarray*}
G_{i\delta}(\beta^u)& = &\sum\limits_{j=1}^{p-1}\beta^{jlp^2}+\sum\limits_{j=0}^{p-1} D_{(i+j)\delta}(\beta)D_{j\delta}(\beta^{lp})\\
            & = &p-1+\sum\limits_{j=0}^{p-1} D_{(i+j)\delta}(\beta) =\sum\limits_{j=0}^{p-1} D_{j}(\beta)=0=s^{(i\delta)}_u.
\end{eqnarray*}
For $u\in  D_{k\delta}$ with $0\le k<p$, we have  by Lemma \ref{lemma1}
 \begin{eqnarray*}
G_{i\delta}(\beta^u)& = &\sum\limits_{j=1}^{p-1}\beta^{jpu}+\sum\limits_{j=0}^{p-1} D_{(i+j)\delta}(\beta)D_{j\delta}(\beta^{u})\\
            & = &1+\sum\limits_{j=0}^{p-1} D_{(i+j)\delta}(\beta)D_{(k+j)\delta}(\beta)\\
            & = &1+C_{i}C_{k}^T\\
            & = &\left\{
\begin{array}{ll}
1, & \mathrm{if}\,\ i=k,\\
0, & \mathrm{otherwise}.
\end{array}
\right. =s^{(i\delta)}_u.
\end{eqnarray*}
That is, $s^{(i\delta)}_u=G_{i\delta}(\beta^u)$ for all $u\ge 0$. We complete the proof. ~\hfill $\square$

\section{Trace Representation}

The trace function from the finite field $\F_{2^n}$ to $\F_{2^k}$ is defined
by
$$
\mathrm{Tr}^n_k(x)=x+x^{2^k}+x^{2^{2k}}+\ldots+x^{2^{(n/k-1)k}}.
$$
For $a,b\in\F_{2^k}$ and $x,y\in\F_{2^n}$, we have $\mathrm{Tr}^n_k(ax+by)=a\mathrm{Tr}^n_k(x)+b\mathrm{Tr}^n_k(y)$. We
refer the reader to \cite{LN} for details on the trace function.

With $\delta$ as above, by Lemma \ref{single-def-pair} it is easy to see that $(G(x),\beta)$ with
\begin{equation}\label{def-pair-e}
G(x)=\sum\limits_{l=\frac{p+1}{2}}^{p-1}G_l(x)=\frac{p-1}{2}\sum\limits_{j=1}^{p-1}x^{jp}+
\sum\limits_{j=0}^{p-1} ~ \sum\limits_{l=\frac{p+1}{2}}^{p-1}D_{l+j\delta}(\beta) D_{j\delta}(x)
\end{equation}
is  the defining pair of $(e_u)$ defined in (\ref{binarythreshold}).

\begin{theorem}\label{trace-e}
Let $q_p(g)=\delta$ for a (fixed) primitive root $g$ modulo $p^2$ and $\beta \in \overline{\mathbb{F}}_{2}$ be a  primitive $p^2$-th root of unity as before. Let $\lambda$ be the smallest positive integer satisfying $2^{\lambda}\equiv 1 \pmod p$. Then the trace representation of  $(e_u)$ defined in
(\ref{binarythreshold}) is
$$
e_u=\frac{p-1}{2}\sum\limits_{k=0}^{\frac{p-1}{\lambda}-1}\mathrm{Tr}^{\lambda}_1(\beta^{upg^k})+\sum\limits_{j=0}^{p-1} \eta_{j}\sum\limits_{k=0}^{\frac{p-1}{\lambda}-1}\mathrm{Tr}^{\lambda p}_p(\beta^{u g^{kp+j}}), ~~ u\ge 0
$$
if $2^{p-1}\not\equiv 1 \pmod {p^2}$, and otherwise
$$
e_u=\frac{p-1}{2}\sum\limits_{k=0}^{\frac{p-1}{\lambda}-1}\mathrm{Tr}^{\lambda}_1(\beta^{upg^k})+\sum\limits_{j=0}^{p-1} \eta_{j}\sum\limits_{k=0}^{\frac{p-1}{\lambda}-1}\mathrm{Tr}^{\lambda}_1(\beta^{u g^{kp+j}}), ~~ u\ge 0,
$$
where
$$
\eta_j=\sum\limits_{l=\frac{p+1}{2}}^{p-1}D_{l+j\delta}(\beta), ~~0\le j<p.
$$
\end{theorem}
Proof. From (\ref{def-pair-e}) we only need to represent $\sum\limits_{j=1}^{p-1}x^{jp}$ and $D_{j\delta}(x) ~(0\le j<p)$ via trace functions. Since  $2^{\lambda}\equiv 1 \pmod p$ and $g$ is also a primitive root modulo $p$, we have
$$
\Z_p^*=\bigcup\limits_{k=0}^{\frac{p-1}{\lambda}-1}g^k \langle 2\rangle,
$$
where $\langle 2\rangle=\{1,2,2^2,\ldots,2^{\lambda-1}\}$ generated by $2$ modulo $p$ is a subgroup of $\Z_p^*$. Hence we derive
$$
\sum\limits_{j=1}^{p-1}x^{jp}=\sum\limits_{k=0}^{\frac{p-1}{\lambda}-1}\mathrm{Tr}^{\lambda}_1(x^{pg^k}).
$$

Now we consider the trace representation of $D_{j\delta}(x) ~(0\le j<p)$. If $2^{p-1}\not\equiv 1 \pmod {p^2}$, we note that the order of $2$ modulo $p^2$ is $\lambda p$ and we find that a power of $2$ modulo $p^2$ belonging to $D_0$ is of the form $2^{ip}~ (0\le i<\lambda)$, all of which are generated by $2^p~(=g^{\frac{p-1}{\lambda}p})$  modulo $p^2$ and form a subgroup of $D_0$. So we have
$$
D_0=\bigcup\limits_{k=0}^{\frac{p-1}{\lambda}-1}g^{kp} \langle 2^p\rangle ~~ \mathrm{and} ~~ D_{j\delta}=g^jD_0=\bigcup\limits_{k=0}^{\frac{p-1}{\lambda}-1}g^{j+kp} \langle 2^p\rangle, ~~0\le j<p,
$$
from which we derive
$$
D_{j\delta}(x)=\sum\limits_{u\in D_{j\delta}}x^u=\sum\limits_{k=0}^{\frac{p-1}{\lambda}-1}\mathrm{Tr}^{\lambda p}_p(x^{g^{kp+j}}), ~~  0\le j<p.
$$
For the case of $2^{p-1}\equiv 1 \pmod {p^2}$, we see that $2\in D_0$ and $\langle 2\rangle=\{1,2,2^2,\ldots,2^{\lambda-1}\}$  generated by $2$ modulo $p^2$ is a subgroup of  $D_0$. One can represent $D_{j\delta}$ as
$$
D_{j\delta}=\bigcup\limits_{k=0}^{\frac{p-1}{\lambda}-1}g^{j+kp} \langle 2\rangle, ~~ 0\le j<p
$$
and hence
$$
D_{j\delta}(x)=\sum\limits_{u\in D_{j\delta}}x^u=\sum\limits_{k=0}^{\frac{p-1}{\lambda}-1}\mathrm{Tr}^{\lambda}_1(x^{g^{kp+j}}), ~~ 0\le j<p.
$$
We complete the proof. ~\hfill $\square$

We need to determine the values of $\eta_j$ for $0\le j<p$. Clearly $\eta_0=E(\beta)$, where $E(x)$ is the generating polynomial of $(e_u)$. For $2^{p-1}\not\equiv 1 \pmod {p^2}$,
it was checked in \cite{CD} that all $\eta_j$'s are different and for each $0\le j<p$ there exists a
$0\le r_j<p$ such that $\eta_j=\eta_0^{2^{r_j}}$. Hence all $\eta_j$'s are (nonzero) roots of a primitive polynomial of degree $p$ over $\F_2$ in this case. Then using the relation (\ref{lc-rho}), we can get an earlier result in \cite{CD} on the linear complexity of $(e_u)$, which we will list it in Section \ref{LC} as a theorem. However, for $2^{p-1}\equiv 1 \pmod {p^2}$, each $\eta_j\in \F_2$ and we have no more results.\\

Similarly, one can obtain  the defining pair $(H(x),\beta)$ of $(f_u)$ defined in (\ref{binarylegendre}).
With $\delta$ as before, from Lemma \ref{single-def-pair} $H(x)$ is of the form
\begin{equation}\label{def-pair-f}
H(x)=\sum\limits_{l\in N}G_l(x)=\frac{p-1}{2}\sum\limits_{j=1}^{p-1}x^{jp}+
\sum\limits_{j=0}^{p-1}~ \sum\limits_{l\in N}D_{l+j\delta}(\beta) D_{j\delta}(x),
\end{equation}
where $N$ is the set of quadratic non-residues modulo $p$.

Now we describe the trace representation of $(f_u)$. The proof is the same as that of Theorem \ref{trace-e} and we omit it.
\begin{theorem}\label{trace-f}
Let $q_p(g)=\delta$ for a (fixed) primitive root $g$ modulo $p^2$ and $\beta \in \overline{\mathbb{F}}_{2}$ be a  primitive $p^2$-th root of unity as before. Let $\lambda$ be the smallest positive integer satisfying $2^{\lambda}\equiv 1 \pmod p$. Then the trace representation of  $(f_u)$ defined in
(\ref{binarylegendre}) is
$$
f_u=\frac{p-1}{2}\sum\limits_{k=0}^{\frac{p-1}{\lambda}-1}\mathrm{Tr}^{\lambda}_1(\beta^{upg^k})+\sum\limits_{j=0}^{p-1} \eta'_{j}\sum\limits_{k=0}^{\frac{p-1}{\lambda}-1}\mathrm{Tr}^{\lambda p}_p(\beta^{u g^{kp+j}}), ~~ u\ge 0
$$
if $2^{p-1}\not\equiv 1 \pmod {p^2}$, and otherwise
$$
f_u=\frac{p-1}{2}\sum\limits_{k=0}^{\frac{p-1}{\lambda}-1}\mathrm{Tr}^{\lambda}_1(\beta^{upg^k})+\sum\limits_{j=0}^{p-1} \eta'_{j}\sum\limits_{k=0}^{\frac{p-1}{\lambda}-1}\mathrm{Tr}^{\lambda}_1(\beta^{u g^{kp+j}}), ~~ u\ge 0,
$$
where for $N$, the set of quadratic non-residues modulo $p$, $\eta'_j$ is
$$
\eta'_j=\sum\limits_{l\in N}D_{l+j\delta}(\beta), ~~0\le j<p.
$$
\end{theorem}

\section{Linear Complexity}\label{LC}

The linear complexity of $(e_u)$ has been shown in \cite{CD}, we present here for completeness of this paper without a proof.

\begin{theorem}\label{lc-e}\cite{CD}
Let $(e_u)$ be the binary sequence of period $p^{2}$ defined in
(\ref{binarythreshold}).  If $2^{p-1}\not\equiv 1 \pmod {p^2}$, then
the linear complexity of $(e_u)$  satisfies
\[
 L((e_u))=\left\{
\begin{array}{ll}
p^2-p, & \mathrm{if}\,\ p \equiv 1 \pmod 4, \\
p^2-1, & \mathrm{if}\,\  p \equiv 3 \pmod 4.\\
\end{array}
\right.\\
\]
\end{theorem}

Fortunately we have a similar result of linear complexity for $(f_u)$. However it seems that the proof in \cite{CD} can't be used and we will give a new simple proof below. We note that the new proof can be used to show Theorem \ref{lc-e}.

\begin{theorem}\label{lc-f}
Let $(f_u)$ be the binary sequence of period $p^{2}$ defined  in
(\ref{binarylegendre}).  If $2^{p-1}\not\equiv 1 \pmod {p^2}$, then
the linear complexity of $(f_u)$  satisfies
\[
 L((f_u))=\left\{
\begin{array}{ll}
p^2-p, & \mathrm{if}\,\ p \equiv 1 \pmod 4, \\
p^2-1, & \mathrm{if}\,\  p \equiv 3 \pmod 4.\\
\end{array}
\right.\\
\]
\end{theorem}
Proof. Let
$$
\Lambda_{\ell}(x)=
\sum\limits_{l\in N}D_{l+\ell}(x)\in
\mathbb{F}_2[x], ~~~\ell=0,\ldots,p-1.
$$
Clearly $\Lambda_0(x)$ is the generating polynomial of $(f_u)$. For $0\le \ell<p$, we only need to show $\Lambda_{\ell}(\beta)\neq 0$  for  a  primitive $p^2$-th root of unity $\beta \in \overline{\mathbb{F}}_{2}$ by  (\ref{def-pair-f}).

We assume that $\Lambda_{\ell_0}(\beta)=0$ for some $0\le \ell_0<p$. Since $2^{p-1}\not\equiv 1 \pmod {p^2}$, we set $q_p(2)=\mu\neq 0$. Then by Lemma \ref{lemma1} we have for $0\le j<p$
$$
0=\left(\Lambda_{\ell_0}(\beta)\right)^{2^j}=\Lambda_{\ell_0}(\beta^{2^j})=\Lambda_{\ell_0+j\mu}(\beta),
$$
where the subscript of $\Lambda$ is reduced modulo $p$. That is, $\Lambda_{\ell}(\beta)=0$ for all $0\le \ell<p$.
We have furtherly $\Lambda_{\ell}(\beta^u)=0$  by Lemma \ref{lemma1} again for all $u\in\Z_{p^2}^*$ and $0\le \ell<p$.
On the other hand, all ($p^2-p$ many) elements $\beta^u$ for $u\in\Z_{p^2}^*$ are roots of
$$
\Phi(x)=1+x^p+x^{2p}+\ldots+x^{(p-1)p}\in \F_2[x],
$$
which has no other roots. Hence  we have
$$
\Phi(x)|\Lambda_0(x) ~~ \mathrm{in}~ \overline{\F}_2[x].
$$
 Let
\begin{equation}\label{pi}
\Lambda_0(x)\equiv \Phi(x)\pi(x) \pmod {x^{p^2}-1}.
\end{equation}
Using the fact that
$$
x^p\Phi(x)  \equiv \Phi(x) \pmod {x^{p^2}-1},
$$
we restrict $\deg(\pi(x))<p$. However, $\Lambda_0(x)$ has $(p-1)^2/2$ terms and the
right hand side of (\ref{pi}) has $pt$ terms if $\pi(x)$ has $t$ terms, a
contradiction.  So we conclude that
$\Lambda_{\ell}(\beta^u)\neq 0$  for all $u\in\Z_{p^2}^*$ and $0\le \ell<p$. Then by (\ref{def-pair-f}) and (\ref{lc-rho})  we get the desired result.

Of course one also can use (\ref{licom}) to finish the proof.
By Lemma \ref{lemma-add}(ii) we get
$$
\Lambda_{\ell}(\beta^{kp})=\left\{
\begin{array}{ll}
0, & \mathrm{if}~ k= 0,\\
(p-1)/2, & \mathrm{if}~ 1\le k<p,
\end{array}
\right.
$$
for $0\le \ell<p$. We draw a conclusion that $\Lambda_0(x)$, the generating polynomial of $(f_u)$, and $x^{p^2}-1$ have exactly $p$ many common roots $\beta^{kp}~(0\le k<p)$ if $p \equiv 1 \pmod 4$ and one common root $\beta^0$ otherwise. ~\hfill $\square$

\section{Concluding Remarks}

In this paper, we explicitly
describe the trace representations of two families of binary
sequences derived from some union of cosets defined using Fermat quotients modulo $p$ via constructing
their defining pairs, which are a sum of the defining pair of binary characteristic sequences defined by a single coset.
We also obtain the linear complexity of these sequences from their defining pairs. In particular, we show that the linear complexity of the Legendre-Fermat sequences $(f_u)$ equals $p^2-p$ or $p^2-1$, depending
 whether $p\equiv 1$ or $3\pmod 4$, under the assumption of $2^{p-1}\not\equiv 1 \pmod {p^2}$. For the case of $2^{p-1}\equiv 1 \pmod {p^2}$, it seems that we can't get more accurate results without additional ideas and we leave it open. But we have mentioned in \cite{CD} that primes $p$ satisfying $2^{p-1}\equiv 1
\pmod {p^2}$ are very rare. To date the only known such primes are
$p=1093$ and $p=3511$ and it was  reported that there are no new
such primes $p< 4\times 10^{12}$, see \cite{CDP1997}.

We  finally remark that, it is more frequent to define balanced binary sequences for some special applications.
So we can modify the definitions of $(e_u)$  and $(f_u)$ to define respectively
$$
\widetilde{e}_u=\left\{
\begin{array}{ll}
0, & \mathrm{if}\,\ u \bmod p^2\in D_0 \cup \cdots \cup D_{(p-1)/2},\\
1, & \mathrm{if}\,\ u \bmod p^2 \in D_{(p+1)/2} \cup \cdots \cup D_{p-1} \cup P,
\end{array}
\right. \quad u\ge 0
$$
and
$$
\widetilde{f}_u=\left\{
\begin{array}{ll}
0, & \mathrm{if}\,\ u\bmod p^2\in \cup_{l\in Q} D_l  \cup D_0,\\
1, & \mathrm{if}\,\ u\bmod p^2\in \cup_{l\in N} D_{l}\cup P,
\end{array}
\right. \quad u\ge 0.
$$
By the fact that
$$
\sum\limits_{j=0}^{p-1}\beta^{jup}=\left\{
\begin{array}{ll}
0, & \mathrm{if}\,\ p\nmid u,\\
1, & \mathrm{otherwise},
\end{array}
\right.
$$
for the given $\beta$, the defining polynomial $\widetilde{G}(x)$ of $(\widetilde{e}_u)$ is of the form
$$
\widetilde{G}(x)=\frac{p-1}{2}\sum\limits_{j=1}^{p-1}x^{jp}+\sum\limits_{j=0}^{p-1}x^{jp}+
\sum\limits_{j=0}^{p-1} ~ \sum\limits_{l=\frac{p+1}{2}}^{p-1}D_{l+j\delta}(\beta) D_{j\delta}(x)
$$
and the defining polynomial $\widetilde{H}(x)$ of $(\widetilde{f}_u)$ is of the form
$$
\widetilde{H}(x)=\frac{p-1}{2}\sum\limits_{j=1}^{p-1}x^{jp}+\sum\limits_{j=0}^{p-1}x^{jp}+
\sum\limits_{j=0}^{p-1}~ \sum\limits_{l\in N}D_{l+j\delta}(\beta) D_{j\delta}(x).
$$
Then it is easy to derive the trace representations and the linear complexity of $(\widetilde{e}_u)$  and $(\widetilde{f}_u)$ from  their defining pairs, respectively.

\section*{Acknowledgements}

The author wishes to thank Ming Su and Arne Winterhof for helpful suggestions.

The work was partially supported  by the National Natural Science
Foundation of China under grant No. 61170246 and Special Scientific Research Program in Fujian Province Universities of China  under grant No. 2013JK044.

Parts of this paper were written during a very pleasant visit of the
 author to RICAM, Austrian Academy of Sciences in Linz. He wishes to thank for the hospitality.

\end{document}